\documentclass[11pt]{amsart}

\usepackage[all]{xy}
\usepackage{amscd}
\usepackage{mathrsfs}
\usepackage{amssymb}
\usepackage{amsthm}
\usepackage{amsmath}

\newtheorem{theorem}{Theorem}[section]
\newtheorem{lemma}[theorem]{Lemma}
\newtheorem{proposition}[theorem]{Proposition}
\newtheorem{corollary}[theorem]{Corollary}
\newtheorem*{ithm}{Theorem}

\theoremstyle{definition}
\newtheorem{definition}[theorem]{Definition}
\newtheorem{example}[theorem]{Example}
\newtheorem{remark}[theorem]{Remark}

\theoremstyle{remark}

\numberwithin{equation}{theorem}

\newcommand{\bde}{\begin{definition}}
\newcommand{\ede}{\end{definition}}
\newcommand{\bpr}{\begin{proposition}}
\newcommand{\epr}{\end{proposition}}
\newcommand{\bth}{\begin{theorem}}
\newcommand{\ethm}{\end{theorem}}
\newcommand{\bexa}{\begin{example}}
\newcommand{\eexa}{\end{example}}
\newcommand{\bcor}{\begin{corollary}}
\newcommand{\ecor}{\end{corollary}}
\newcommand{\blem}{\begin{lemma}}
\newcommand{\elem}{\end{lemma}}
\newcommand{\brem}{\begin{remark}}
\newcommand{\erem}{\end{remark}}
\newcommand{\bprf}{\begin{proof}}
\newcommand{\eprf}{\end{proof}}
\def \benu{\begin{enumerate}\renewcommand{\labelenumi}{(\roman{enumi})}\renewcommand{\itemsep}{0pt}}
\newcommand{\eenu}{\end{enumerate}}
\newcommand{\fkp}{{\mathfrak p}}

\newcommand{\mE}{{\mathcal E}}
\newcommand{\mF}{{\mathcal F}}

\newcommand{\mI}{{\mathcal I}}

\newcommand{\mM}{{\mathcal M}}

\newcommand{\mO}{{\mathcal O}}

\newcommand{\N}{\mathbb{N}}

\newcommand{\Z}{\mathbb{Z}}

\newcommand{\Q}{\mathbb{Q}}
\newcommand{\F}{\mathbb{F}}

\newcommand{\ke}{{\rm Ker\,}}

\def\Gal{{\rm Gal}}

\def\Pic{{\rm Pic}}

\DeclareMathOperator{\Spec}{\mathrm{Spec}}

\newcommand{\Div}{{\rm Div}}

\def\div{{\rm div}}

\newcommand{\bsc}{\begin{tobira}}
\newcommand{\esc}{\end{tobira}}
\newcommand{\abc}{\xymatrix{0 \ar@{>}[r] & A \ar@{>}[r]^{f} & B \ar@{>}[r]^{g} & C \ar@{>}[r] & 0}}
\newcommand{\xyz}{\xymatrix{0 \ar@{>}[r] & X \ar@{>}[r]^{f} & Y \ar@{>}[r]^{g} & Z \ar@{>}[r] & 0}}
\def\phi{\varphi}

\def\lra{\longrightarrow}
\def\Lra{\Longrightarrow}

\def\ra{\rightarrow}

\def\ol{\overline}

\def\ul{\underline}

\DeclareMathOperator{\chara}{\mathrm{char}}

\newcommand{\iso}{\stackrel{\sim}{\ra}}
\newcommand{\liso}{\stackrel{\sim}{\lra}}

\newcommand{\pr}{{\operatorname{pr}}}

\begin{document}

\title[Correspondences on weight spectral sequences]{On the action of algebraic correspondences on weight spectral sequences}

\author[T. Yoshida]{Teruyoshi Yoshida}
\address{University of Cambridge, Department of Pure Mathematics and Mathematical Statistics, Centre for Mathematical Sciences, Wilberforce Road, Cambridge, CB3 0WB, UK}
\email{t.yoshida\char`@dpmms.cam.ac.uk}
\date{\today}

\begin{abstract}
In a work of T.~Saito, the action of algebraic correspondences on the \'etale cohomology of varieties over local fields with semistable reduction is related to correspondences on smaller strata via weight spectral sequences. We give an intersection theoretic construction of these correspondences. Under a finiteness condition this enables us to compute them without involving the blow-ups of products, and prove their compatibility with compositions. These features are essential for the application to Shimura varieties.
\end{abstract}

\maketitle

\tableofcontents

\section{Introduction}

In this paper we study the construction of T.~Saito \cite{SaT} on the action of algebraic correspondences on the weight spectral sequences, which compute the $\ell$-adic \'etale cohomology of varieties over local fields with semistable reduction. Using a blow-up of a product of semistable schemes, Saito constructs algebraic correspondences on smaller strata which act compatibly on the terms of weight spectral sequences. We prove that these correspondences can be obtained by intersection theory on the integral model, and that they satisfy an intersection theoretic formula which does not involve the blow-ups. Under a finiteness condition this formula enables us to compute them directly, and to prove their compatibility with compositions. In \cite{Y}, we will apply this result to Hecke correspondences on certain Shimura varieties.

Let $K$ be a complete discrete valuation field with a finite residue field $k$ and the ring of integers $\mO$, and $X$ a strictly semistable scheme (see \S\ref{sectionsst} for the definition) of relative dimension $n-1$ over $\mO$. Then its special fiber $Y:=X\times_{\mO}k$ is written as $Y=\bigcup_{i\in \Delta}Y_i$ with $\Delta:=\{1,...,t\}$ and $Y_i$ proper smooth over $k$, where $Y_i$ and $Y_j$ intersect transversally for $i\neq j$. Let $Y_I:=\bigcap_{i\in I}Y_i$ for any finite subset $I\subset \Delta$, and $Y^{(m)}:=\coprod_{|I|=m}Y_I$ for $1\leq m\leq n$. For a prime $\ell\neq \chara k$, the {\em weight spectral sequence} (\cite{RZ},\cite{SaT}) reads
\[E_1^{i,j}:=\bigoplus_{s\geq \max(0,-i)}H^{j-2s}\bigl(Y^{(i+2s+1)}_{\ol{k}},\ \Q_\ell(-s)\bigr) \Lra H^{i+j}(X_{\ol{K}}, \Q_\ell).\]

%Now let $\Gamma$ be an algebraic correspondence on $X$, namely an $n$-dimensional cycle on $X\times_{\mO}X$ which is flat over $\mO$. 

Let $\Gamma_K$ be an algebraic corerspondence on $X_K$, namely an $(n-1)$-dimensional cycle on $X_K\times X_K$. We are interested in its action on $H^*(X_{\ol{K}},\Q_\ell)$, defined as $[\Gamma_K]^*:=\pr_{1*}\circ ([\Gamma_K]\cup)\circ \pr_2^*$. We denote $Y_{I,J}:=Y_I\times_k Y_J$ for $I,J\subset \Delta$, and in particular $Y_{i,j}:=Y_{\{i\},\{j\}}$. Let $(X\times_{\mO}X)_{\rm sm}$ be the smooth locus of the morphism $X\times_{\mO}X\lra \Spec \mO$, and let $Y_{i,j}^0:=Y_{i,j}\cap (X\times_{\mO}X)_{\rm sm}$. Then $Y_{i,j}^0$ is a Cartier divisor of $(X\times_{\mO}X)_{\rm sm}$. Let $\Gamma$ be the closure of $\Gamma_K$ in $X\times_\mO X$. Our main result is stated as follows:

\begin{ithm}
Assume that two projection maps $\Gamma\ra X$ are both finite morphisms.
\begin{enumerate}
\renewcommand{\labelenumi}{(\roman{enumi})}
\item There is a unique collection $\{\Gamma_{I,J}\}$ of cycles $\Gamma_{I,J}$ on $Y_{I,J}$ for all pairs $(I,J)$ with $|I|=|J|$, satisfying the following two conditions.
\begin{itemize}
\item $\Gamma_{i,j}$ is the closure of the cycle $\Gamma_{i,j}^0:=Y_{i,j}^0\cdot \Gamma|_{(X\times_\mO X)_{\rm sm}}$ in $Y_{i,j}$.
\item When $|I|=|J|+1=m$ and $I=\{i_1,...,i_m\},\ J=\{j_1,\ldots,j_{m-1}\}$ are in increasing order, there is an equality:
\[ \sum_{h=1}^m(-1)^h\,Y_{I,J}\cdot \Gamma_{I\setminus \{i_h\},J} = \sum_{j\in \Delta\setminus J} (-1)^{h(j)}\,\Gamma_{I,J\cup\{j\}} \]
of $(n-m)$-dimensional cycles on $Y_{I,J}$, where $1\leq h(j)\leq m$ is determined by $j_{h(j)-1}<j<j_{h(j)}$ (set $j_m:=\infty$). 
\end{itemize}
Then setting $\Gamma^{(m)}:=\coprod_{|I|=|J|=m}\Gamma_{I,J}$ as an $(n-m)$-dimensional cycle on $Y^{(m)}\times_k Y^{(m)}$ for $1\leq m\leq n$, the action $\oplus [\Gamma^{(i+2s+1)}]^*$ on $E_1^{i,j}$ is compatible with the action $[\Gamma_K]^*$ on $H^{i+j}(X_{\ol{K}},\Q_\ell)$ via the weight spectral sequence.
\item If $\Gamma_1,\Gamma_2$ are two correspondences as above, then $(\Gamma_2\circ \Gamma_1)^{(m)}=\Gamma_2^{(m)}\circ \Gamma_1^{(m)}$ for $1\leq m\leq n$ (here $\Gamma_2\circ \Gamma_1$ is the closure of $\Gamma_{2,K}\circ \Gamma_{1,K}$ in $X\times_\mO X$).
\eenu
\end{ithm}

In particular, in a situation where we have an algebra consisting of correspondences (e.g.\ a Hecke algebra consisting of Hecke correspondences on Shimura varieties), the second part of the theorem shows that the algebra acts on each term of the weight spectral sequence and the sequence is equivariant. For Hecke correspondences on Shimura varieties the finiteness condition is rather restrictive, but it is satisfied in the cases treated in \cite{TY},\cite{Y}.

The cycles $\Gamma_{I,J}$ in the part (i) of the theorem are similar to the ones defined in \cite{SaT} using a blow-up of $X\times_\mO X$. We first prove that these cycle classes (up to rational equivalence) are obtained by intersection theory (Definition \ref{cycleclasses}, Proposition \ref{wssaction}). This was also proved in \cite{Mi}, \S3.2 by a slightly different method; it improves the construction in \cite{SaT} by removing a denominator. Then we prove the formula in the theorem by a repeated use of the projection formula, hence it is valid for general $\Gamma$ as an equality of cycle classes (Proposition \ref{iteratedproj}). Then we impose the finiteness condition in order to ensure that this formula characterizes (in fact computes) the cycles $\Gamma_{I,J}$, which makes them rigid enough to satisfy (ii).

\subsection*{Acknowledgment}

This work is based on a part of the PhD thesis of the author submitted to Harvard University in 2006. The author would like to thank his adviser Richard Taylor for his encouragements; computations in the curve case is due to him. The question of computing the correspondences $\Gamma^{(m)}$ was originally raised for the purpose of proving the results in \cite{TY}. The author thanks Yoichi Mieda and Seidai Yasuda for helpful discussions.

\subsection*{Notation}

For a field $F$ we denote its separable closure by $\ol{F}$. We often suppress the notation $\Spec$ from the fiber products of/over affine schemes, and denote a base change to a (geometric) point of the base scheme by a subscript (e.g.\ $X_{\ol{K}}=X\times_K\ol{K}$ means $X\times_{\Spec K}\Spec \ol{K}$). For closed subschemes $V,W$ of a scheme $X$, we write $V\cap W:=V\times_XW$. If $V,W$ are reduced, then $V\cup W$ is the reduced closed subscheme with the underlying space $V\cup W$. A scheme is {\em integral} if it is irreducible and reduced. For a scheme $X$, its associated reduced subscheme is denoted by $X_{\rm red}$.

\section{Review of some intersection theory} \label{review}

We include a summary of what we need from basic intersection theory on noetherian schemes, and fix the detailed notation. Our reference is \cite{Fu2}, Chapters 1 and 2, where the schemes are assumed to be over a field, but as is remarked in Chapter 20 of \cite{Fu2}, what we recall here are proven for noetherian schemes (see also \cite{Fu1}, \S1). Up to \S\ref{charcorres} we only need intersections of cycles with divisors, but in \S\ref{composesec} we will use the (refined) intersection products of cycles (\cite{Fu2}, Chapters 6/8 and \S20.2), mainly over a field.

Let $X,Y,$ etc.\ be noetherian separated equidimensional schemes, assumed to be over a regular base ring so that all rings are universally catenary.

\subsection{Cycles} \label{sscycles}

For $i\in \N$, an $i$-dimensional {\it cycle} on $X$ is an element of the free abelian group $Z_i(X)$ generated by the elements $[V]$ for all integral closed subschemes $V$ of dimension $i$ in $X$. We have $Z_i(X)\cong Z_i(X_{\rm red})$ for all $i$. Let $Z(X):=\bigoplus_{i\geq 0}Z_i(X)$.

For a coherent sheaf $\mF$ on $X$, let $W_1,...,W_m$ be the irreducible components of the support of $\mF$. If $x_j$ is the generic point of $W_j$ and $\mO_{x_j}$ is the local ring of $X$ at $x_j$, then {\it the cycle associated to $\mF$} is defined by
\[Z(\mF):=\sum_{j=1}^m{\rm length}_{\mO_{x_j}}(\mF_{x_j})\cdot [W_j]\in Z(X).\]
When $\dim W_j\leq i$ for all $j$, then set $Z_i(\mF)\in Z_i(X)$ to be the sum of terms for $W_j$ with $\dim W_j=i$. For any closed subscheme $Y$ of $X$, set $[Y]:=Z(\mO_Y)$.%, and when $\dim Y\leq i$ then set $[Y]_i:=Z_i(\mO_Y)$. 

If $\Gamma=\sum_jm_jV_j\in Z(X)$ where $m_j\in \Z\setminus\{0\}$ and $V_j$ are integral closed subschemes of $X$, then we call $|\Gamma|:=\bigcup_jV_j$ the {\it support} of $\Gamma$. It is a reduced closed subscheme of $X$; e.g.\ for a closed subscheme $Y$ of $X$ we have $|[Y]|=Y_{\rm red}$. 

For a {\em proper} morphism $f:X\ra Y$ between noetherian schemes, the {\it push-forward} $f_*:Z_i(X)\ra Z_i(Y)$ is defined.
In particular, if $V$ is a closed subscheme of $X$, we often identify $Z_i(V)$ with a subgroup of $Z_i(X)$ via push-forward map. 

For a {\em flat} morphism $f:X\ra Y$ between noetherian schemes, the {\it pull-back} $f^*:Z(Y)\ra Z(X)$ is defined by $f^*[V]:=[V\times_YX]$ for an integral closed subscheme $V$ of $Y$, and extending linearly. We have $f^*[V]=[V\times_YX]$ for any closed subscheme $V$ of $Y$. 

Let $j:U\ra X$ be an open immersion. We write $\Gamma|_U:=j^*(\Gamma)$ for $\Gamma\in Z(X)$. For a proper morphism $f:X\ra Y$, an open subscheme $U\subset Y$ and $\Gamma\in Z(X)$, we have $f_*(\Gamma)|_U=f'_*(\Gamma|_{f^{-1}(U)})$ where $f':=f|_{f^{-1}(U)}$. For a cycle on $U$, we define its {\it closure} in $X$ by linearly extending the scheme-theoretic closure operation on integral subschemes. This gives a section of $j^*$, hence $j^*$ is surjective. Thus we have:

\blem \label{opencycles}
If $V$ is the reduced closed subscheme of $X$ with the underlying space $X \setminus U$ and $\iota: V\ra X$ is the closed immersion, then the following is exact:
\[\xymatrix{
0 \ar[r] & Z(X\setminus U)\ar[r]^-{\iota_*} & Z(X) \ar[r]^-{j^*} & Z(U) \ar[r] & 0
}\]
In particular, if $\dim (X\setminus U)<i$, then $j^*$ is injective on $Z_i(X)$.
\elem

\subsection{Divisors and cycle classes}

We write $\Div(X):=H^0(X,\mM_X^\times/\mO_X^\times)$, the group of {\it Cartier divisors} on $X$, where $\mM_X$ is the sheaf of total quotient rings of $\mO_X$ and $^\times$ means the sheaf of groups consisting of invertible elements. We write the multiplication of sections additively. For $D\in \Div(X)$, the {\em support} $|D|$ of $D$ is the reduced closed subscheme of $X$ where $D$ is not a section of $\mO_X^\times$. A Cartier divisor is {\it principal} if it comes from a global section of $\mM_X^\times$, and let $\Pic(X)$ be the group of Cartier divisors modulo principal ones. We denote by $\mO(D)$ the invertible sheaf (up to isomorphism) which corresponds to the image of $D$ under the coboundary map $H^0(X,\mM_X^\times/\mO_X^\times)\lra H^1(X,\mO_X^\times)$. Thinking of $D$ as a collection of local sections of $\mM_X$, the sheaf $\mO(D)$ is realized as a $\mO_X$-submodule of $\mM_X$ locally generated by $-D$. On the open subscheme $U:=X\setminus |D|$, the section $1\in H^0(U,\mM_X)$ gives a section $s_D\in H^0(U,\mO(D))$, or a trivialization $s_D:\mO_U\cong \mO(D)|_U$. We can recover $D$ from the pair $(\mO(D),s_D)$ (\cite{Fu2}, Lemma 2.2(a)). 

If $D$ is {\it effective}, i.e.\ if $D$ comes from a collection of local sections of $\mO_X$, or equivalently $\mO_X\subset \mO(D)$, then the section $s_D\in H^0(U,\mO(D))$ extends to a global section $s_D\in H^0(X,\mO(D))$.
% and $|D|$ is defined by $s_D=0$. 
In this case $s_D$ gives an identification of $\mO(-D)$ with a sheaf of ideals of $\mO_X$, and we also denote the corresponding closed subscheme of $X$ by $D$, so that $\mO_D\cong \mO_X/\mO(-D)$. If $\dim X=n$ then we define the associated cycle by $[D]:=Z_{n-1}(\mO_D)\in Z_{n-1}(|D|)$; as $\dim D=n-1$ (\cite{EGA} IV, 21.2.12), it coincides with the class $[D]$ where $D$ is considered as a closed subscheme of $X$, and $|D|=|[D]|=D_{\rm red}$. By linearly extending the map $D\mapsto [D]$ we have a homomorphism $\Div(X)\ra Z_{n-1}(X)$ (\cite{EGA} IV, 21.6.7). It is injective if $X$ is normal (\cite{EGA} IV, 21.6.9). 

For a morphism $f:X\ra Y$ of noetherian schemes and a Cartier divisor $D$ on $Y$, the {\em pull-back} Cartier divisor $f^*D$ on $X$ is defined if $f$ is flat or surjective (\cite{EGA} IV, 21.4.5). If $j:U\ra X$ is an open immersion, we write $D|_U:=j^*D$. 
 
For an $(i+1)$-dimensional integral scheme $V$ and a rational function $\phi\in k(V)^\times=H^0(V,\mM_V^\times)$, we have the associated Cartier divisor $\div(\phi)$ and $[\div(\phi)]\in Z_i(V)$. The group of $i$-dimensional {\it cycle classes} $A_i(X)$ is defined as the cokernel of $\bigoplus_Vk(V)^\times \stackrel{\div}{\lra} Z_i(X)$, where $V$ runs through all $(i+1)$-dimensional integral subschemes of $X$. Let $A(X):=\bigoplus_{i\geq 0}A_i(X)$. If $\dim X=n$ then $Z_n(X)=A_n(X)$, and the map $D\mapsto [D]$ induces a homomorphism $\Pic(X)\ra A_{n-1}(X)$. We have $A_i(X)\cong A_i(X_{\rm red})$ for all $i$. 

For a proper (resp.\ flat) morphism $f:X\ra Y$, the push-forward (resp.\ pull-back) of cycles induces $f_*:A_i(X)\ra A_i(Y)$ (resp.\ $f^*:A(Y)\ra A(X)$). We sometimes denote a cycle class and its image under a push-forward via closed immersions by the same symbol.

\subsection{Intersecting with divisors} \label{intersectdiv}

For a Cartier divisor $D$ on $X$ and an $i$-dimensional integral closed subscheme $V$ of $X$ with the inclusion $\iota:V\ra X$, we will define 
\[D\cdot [V]:=[\iota^*D]\in A_{i-1}(|D|\cap V).\]
Let $\mO(\iota^*D):=\mO(D)\otimes_{\mO_X}\mO_V$ be the invertible sheaf on $V$. There are two possibilities: (1) if $V\subset |D|$, then $V=|D|\cap V$ is $i$-dimensional and $[\iota^*D]$ is defined as the image of $[\mO(\iota^*D)]\in \Pic(V)$ in $A_{i-1}(V)$; (2) if $V$ is not contained in $|D|$, then we define $\iota^*D$ to be the Cartier divisor on $V$ corresponding to $(\mO(\iota^*D), \iota^*(s_D))$, where $\iota^*(s_D)$ is the trivialization $\mO_U\cong \mO(\iota^*D)|_U$, where $U:=V\setminus |D|$, obtained by restricting $s_D$. Hence $|\iota^*D|=(|D|\cap V)_{\rm red}$, and in this case $|D|\cap V$ is $(i-1)$-dimensional and $[\iota^*D]\in Z_{i-1}(|D|\cap V)=A_{i-1}(|D|\cap V)$. For any $\Gamma\in Z_i(X)$, we define $D\cdot \Gamma \in A_{i-1}(|D|\cap |\Gamma|)$ by extending the above definition linearly.

We analyze the situation where $D$ is effective a little further.

\blem \label{tor}
Let $D$ be an effective Cartier divisor $D$ on $X$ and $V$ is an $i$-dimensional integral closed subscheme of $X$ with the inclusion $\iota:V\ra X$.
\benu
\item $\ul{\rm Tor}_i(\mO_D,\mO_V)=0$ for $i>1$ (in this context all $\ul{\rm Tor}$ and $\otimes$ will be over $\mO_X$).
\item If $V$ is not contained in $|D|$, then $\ul{\rm Tor}_1(\mO_D,\mO_V)=0$, and $\iota^*D$ is an effective Cartier divisor on $V$ whose corresponding subscheme is $D\cap V$.
\eenu
\elem

\bprf
(i) Take the long exact sequence for $\ul{\rm Tor}_i(-,\mO_V)$ from the short exact sequence:
\[\xymatrix{
0 \ar[r] & \mO(-D) \ar[r] & \mO_X \ar[r] & \mO_D \ar[r] & 0
}\]
and note that $\mO(-D)$ and $\mO_X$ are locally free, hence flat and have trivial $\ul{\rm Tor}_i$ for $i>0$. (ii) The long exact sequence considered in (i) reads:
\begin{equation} \label{torseq}
\xymatrix{
0 \ar[r] & \ul{\rm Tor}_1(\mO_D,\mO_V) \ar[r] & \mO(-D)\otimes \mO_V \ar[r] & \mO_V \ar[r] & \mO_D\otimes \mO_V \ar[r] & 0.
}\end{equation}
For $\ul{\rm Tor}_1$ the question is local, so we can assume $X=\Spec(A),\ V=\Spec(A/\fkp)$ and $D=\Spec(A/(a))$ with $a\notin \fkp$. Then $\ul{\rm Tor}_1(\mO_D,\mO_V)=\ke \bigl((a)/(a)\fkp \ra A/\fkp\bigr)=\bigl((a)\cap \fkp\bigr)/(a)\fkp$, but
%the above exact sequence is:
%\[\xymatrix{
%0 \ar[r] & (a)\cap \fkp/(a)\fkp \ar[r] & (a)/(a)\fkp \ar[r] & A/\fkp \ar[r] & A/((a)+\fkp) \ar[r] & 0,
%}\]
as $a\notin \fkp$ and $\fkp$ is prime $(a)\cap \fkp=(a)\fkp$. Now note that $\iota^*D$ is a Cartier divisor corresponding to $(\mO(\iota^*D),\iota^*(s_D))$, and $\iota^*(s_D)$ extends to a global section of $\mO(\iota^*D)=\mO(D)\otimes \mO_V$. The map $\mO(-\iota^*D)=\mO(-D)\otimes \mO_V\ra \mO_V$, which was just shown to be injective, corresponds to the trivialization $\iota^*(s_D)$. Hence $\iota^*D$ is an effective Cartier divisor on $V$ with the corresponding subscheme $D\cap V$, as $\mO_{D\cap V}=\mO_D\otimes \mO_V\cong \mO_V/\mO(-\iota^*D)$ (see also \cite{Fu2}, Example 2.6.5 and Lemma A.2.7, and \cite{Fu1}, \S1.7).\eprf

\bpr  \label{projformula}
{\em (\cite{Fu2}, Prop.\ 2.3 and Cor.\ 2.4.2)} If $D,D'$ are Cartier divisors on $X$, then:
\benu
\item $\Gamma \longmapsto D\cdot \Gamma$ induces $A_i(Y)\ra A_{i-1}(|D|\cap Y)$ for any closed subscheme $Y\subset X$.
\item $D \longmapsto D\cdot \Gamma$ induces $\Pic(X)\ra A_{i-1}(|\Gamma|)$, denoted by $[D]\longmapsto [D]\cdot \Gamma$.
\item $D\cdot (\Gamma+\Gamma')=D\cdot \Gamma + D\cdot \Gamma'$ in $A_{i-1}(|D|\cap(|\Gamma|\cup|\Gamma'|))$ for $\Gamma,\Gamma'\in Z_i(X)$.
\item $(D+D')\cdot \Gamma = D\cdot \Gamma + D'\cdot \Gamma$ in $A_{i-1}((|D|\cup |D'|)\cap |\Gamma|)$ for $\Gamma\in Z_i(X)$.
\item $D\cdot (D'\cdot \Gamma)=D'\cdot (D\cdot \Gamma)$ in $A_{i-2}(|D|\cap |D'|\cap |\Gamma|)$ for $\Gamma\in Z_i(X)$.
\item ({\it projection formula}) Let $f:X'\ra X$ be a proper surjective morphism and $\Gamma\in A_i(X')$. Write $g:=f|_{f^{-1}(|D|)\cap |\Gamma|}:f^{-1}(|D|)\cap |\Gamma|\lra |D|\cap f(|\Gamma|)$. Then:
\[g_*(f^*D\cdot \Gamma)=D\cdot f_*(\Gamma)\ \ \ \text{ in } A_{i-1}(|D|\cap f(|\Gamma|)).\]
\item If $U$ is an open subscheme of $X$ and $\Gamma\in Z(X)$, then $D|_U\cdot \Gamma|_U=(D\cdot \Gamma)|_{U'}$, where $U':=U\cap (|D|\cap |\Gamma|)$.
\eenu
\epr

\brem
The formula (vi) makes sense for any proper $f$ by noting that if $\iota:V\ra X$ is any closed immersion we have $[\iota^*D]\cdot \Gamma = D\cdot \iota_*(\Gamma)$ for $\Gamma\in Z(V)$, where $[\iota^*D]:=D\cdot [V]$.
\erem

\subsection{Schemes over Dedekind rings}

Now we consider cycles on schemes over a connected {\it regular 1-dimensional} noetherian base scheme $S$. Let $\eta$ be the generic point of $S$, i.e.\ $S$ is the closure of $\eta$. Let $X$ be a noetherian scheme over $S$ and $\Gamma\in Z_i(X)$. Then $j:X_\eta:=X\times_S\eta\lra X$ is an open immersion, hence we define
\[\Gamma_{\eta}:=j^*(\Gamma)\in Z_i(X_\eta).\] 
For a cycle on $X_\eta:=X\times_S\eta$, taking its closure in $X$ will increase its dimension by one. Observe that for an integral closed subscheme $V$ of $X$, it is the scheme theoretic closure of $V_\eta:=V\times_S\eta$ if and only if $V_\eta\neq \emptyset$, which is equivalent to $V$ being flat over $S$. We call $\Gamma\in Z_i(X)$ a {\it flat} cycle of $X$ if $|\Gamma|$ is flat over $S$, which is equivalent to say that $\Gamma$ is the closure of $\Gamma_\eta$. If $s$ is a closed point of $S$, then $X_s:=X\times_Ss$ is an effective Cartier divisor on $X$, hence we define
\[\Gamma_s:=X_s\cdot \Gamma\in A_{i-1}(X_s\cap|\Gamma|).\]
If $\Gamma\in Z_i(X)$ is a {\em flat} cycle, then $|\Gamma|$ is not contained in $(X_s)_{\rm red}$ (the case (2) of \S\ref{intersectdiv}), hence we have $\Gamma_s\in Z_{i-1}(X_s\cap|\Gamma|)$.

\subsection{Algebraic correspondences} \label{algcor}

In this paper, an {\it algebraic correspondence} on a noetherian $S$-scheme $X$ (or from $X$ to $Y$), where $S$ is the base scheme, is a cycle on the product $X\times_S X$ (resp.\ on $X\times_SY)$ with the dimension equal to $\dim X$. For example, the {\it graph} $\Gamma_f$ associated to an $S$-morphism $f:X\ra Y$ is an algebraic correpondence.

An algebraic correspondence induces a map on the $\ell$-adic \'etale cohomology groups, where we assume $\ell$ is invertible on $S$, as follows. Assume that $S$ is a spectrum of a field or a regular 1-dimensional noetherian scheme, and fix a geometric point $\iota:s\ra S$, where $s$ is a spectrum of a separable closure of the residue field of $\iota(s)$. If $X$ is proper and smooth over $S$ and $Y$ has relative dimension $n$ over $S$, and $\Gamma$ is a correspondence from $X$ to $Y$, we define $[\Gamma]^*$ to be the composite (this coincides with $f^*$ when $\Gamma=\Gamma_f$ for a proper $f$):
\[\xymatrix{
H^i(Y_s, \Q_\ell) \ar[r]^-{\pr_2^*} & H^i(X_s\times Y_s,\ \Q_\ell) \ar[r]^-{{\rm cl}(\Gamma_{\iota(s)})\cup} & H^{i+2n}(X_s\times Y_s,\ \Q_\ell)(n) \ar[r]^-{\pr_{1*}} & H^i(X_s, \Q_\ell),
}\]
where $\Gamma_{\iota(s)}:=\Gamma\times \iota(s)$ and ${\rm cl}$ is the cycle class map (see the proof of Proposition \ref{wssaction} for one construction). 
%When $\Gamma$ is a graph $\Gamma_f$ of a proper $S$-morphism $f:X\ra Y$, the map $[\Gamma_f]^*$ coincides with the pull-back $f^*$ by $f$. 
As the geometric constructions are done over $\iota(s)$, all the maps between cohomology groups commute with the action of $\Gal(s/\iota(s))$.

\section{Correspondences on semistable schemes} \label{sectionsst}

Let $K$ be a complete discrete valuation field with a finite residue field $k\cong \F_q$ with $\chara k=p>0$. Let $\mO$ be the ring of integers of $K$, and $\varpi$ be a uniformizer of $\mO$. Let $X$ be a {\it strictly semistable scheme} over $\mO$ of relative dimension $n-1$, i.e.\ (1) $X$ is proper over $\mO$ and $X_K:=X\times_\mO K$ is smooth over $K$, (2) $X$ is Zariski locally \'etale over $\Spec \mO[X_1,...,X_n]/(\varpi-X_1\cdots X_m)$ for some $1\leq m\leq n$. This implies that the special fiber $Y:=X\times_\mO k$ can be written as a union $Y=\bigcup_{i\in \Delta}Y_i$ where $\Delta=\{1,...,t\}$ for some positive integer $t$, each $Y_i$ is proper smooth of dimension $n-1$ over $k$, and for $i \neq j$ the schemes $Y_i$ and $Y_j$ share no common connected component.

Now let $\Gamma_K$ be an algebraic correspondence on $X_K$, i.e.\ $\Gamma_K\in Z_{n-1}(X_K\times_KX_K)$. We are interested in the action of $\Gamma_K$ on the $\ell$-adic \'etale cohomology of $X_{\ol{K}}$, namely
\[[\Gamma_K]^*:=\pr_{1*}\circ ([\Gamma_K]\cup)\circ \pr_2^*\ \ \text{ on }\ \ H^j(X_{\ol{K}},\Q_\ell),\]
which we investigate via the weight spectral sequence. In \cite{SaT}, T.~Saito constructed a strictly semistable resolution of the self-product 
\[f:X'\lra X\times_\mO X,\]
which is an isomorphism on the generic fibers, and used the closure $\Gamma'$ of $\Gamma_K$ in $X'$ to define correspondences on smaller strata. We will recover his construction via intersection theory, and derive a formula among them which does not refer to the blow-up $X'$. 

\subsection{Semistable resolution of the self-product of semistable schemes}

First we recall the construction of a resolution $f:X'\ra X\times_\mO X$ from \S1.2 of \cite{SaT}. Let $Y:=X\times_\mO k=\bigcup_{i\in \Delta}Y_i$ with $Y_i$ proper smooth of dimension $n-1$ over $k$. For $i\in \Delta$, let $\mI_i$ be the sheaf of ideals of $\mO_X$ corresponding to the closed subscheme $Y_i$ of $X$. Note that each $Y_i$ is a Cartier divisor of $X$ with $\mI_i=\mO(-Y_i)$. For a subset $I\subset \Delta$, let $Y_I:=\bigcap_{i\in I}Y_i$, which is proper smooth of dimension $n-|I|$ over $k$ if not empty. The smooth locus $X_{\rm sm}$ of $X\ra \Spec \mO$ is the complement of $X_{\rm sing}:=\bigcup_{|I|=2}Y_I$ in $X$.

Let $Y_{i,j}:=Y_i\times_kY_j$ for $i,j\in \Delta$. This is a $2(n-1)$-dimensional proper smooth subvariety of $(X\times_\mO X)\times_\mO k$, which in turn is the union of $Y_{i,j}$ for all $(i,j)\in \Delta\times \Delta$. Note that in general $Y_{i,j}$ is not a Cartier divisor of $X\times_\mO X$. We define a partial order on $\Delta\times \Delta$ as follows:
\[ (i,j)\leq (i',j') \iff i\leq i' \text{ and } j\leq j'.\]

\bpr \label{basicblowup}
{\em (Lemma 1.9 of \cite{SaT})} We define a proper morphism $f:X'\ra X\times_\mO X$ as the blow-up of $X\times_\mO X$ by the ideal:
\[\prod_{(i,j)\in \Delta\times \Delta}\Bigl(\prod_{h=1}^i\pr_1^*\mI_h+\prod_{h=1}^{j}\pr_2^*\mI_{h}\Bigr).\]
\benu
\item  Then $X'$ is a strictly semistable scheme over $\mO$, and $f$ is an isomorphism outside $Z:=X_{\rm sing}\times_kX_{\rm sing}$. 
\item Moreover, let $D_{i,j}$ be the closure in $X'$ of $Y_{i,j}\setminus Z$. Then $D_{i,j}$ is proper smooth over $k$, and the special fiber $Y':=X'\times_\mO k$ is the union $Y'=\bigcup_{(i,j)\in \Delta\times \Delta}D_{i,j}$.
\item Let $(i,j),(i',j')\in \Delta\times \Delta$. If $D_{i,j}\cap D_{i',j'}\neq \emptyset$, then either $(i,j)\leq(i',j')$ or $(i,j)\geq (i',j')$. 
\eenu
\epr

\blem \label{cartdivsum}
For all $i\in \Delta$, we have $f^*(\pr_1^*Y_i)=\sum_{j\in \Delta}D_{i,j}$ as Cartier divisors on $X'$.
\elem

\bprf
As $X'$ is regular, hence normal, it is enough to show the equality as $(2n-2)$-dimensional cycles on $X'$. These cycles are equal when restricted to the open subscheme $X'\setminus Z$, and $Z_{2n-2}(X')\iso Z_{2n-2}(X'\setminus Z)$ by Lemma \ref{opencycles}, because $Z$ is at most $(2n-4)$-dimensional.
\eprf

%\bpr
%Let $|I|=|I'|=j$. The proper morphism $f|_{D_{I,I'}}:D_{I,I'}\lra Y_{I,I'}$ is a $(\Ps^1)^{j-1}$-bundle over $Y_{I,I'}$ outside the inverse image of $\bigcup Y_{J,J'}$, where the $J,J'$ runs through all subsets of $\Delta$ with $|J|=|J'|=j+1$.
%\epr

\subsection{Cycle classes on the strata given by a correspondence}

Now consider $\Gamma_K\in Z_{n-1}(X_K\times_KX_K)$ and let $\Gamma\in Z_n(X\times_\mO X),\ \Gamma'\in Z_n(X')$ be the closures of $\Gamma_K$ in $X\times_\mO X$ and $X'$ respectively. First note the following:

\blem
Let $S,S'$ be schemes over $\mO$, and $f:S'\ra S$ a proper morphism which is an isomorphism on the generic fibers. Let $\Gamma_K$ be a cycle on $S_K:=S\times_\mO K$, and $\Gamma,\Gamma'$ be the closures of $\Gamma_K$ in $S,S'$ respectively. Then $|\Gamma|=f(|\Gamma'|)$ and $\Gamma=f_*(\Gamma')$.
\elem

\bprf
It is enough to show that for an integral closed subscheme $V_K$ of $S_K$, if we denote the closures of $V_K$ in $S,S'$ respectively by $V,V'$, then $V=f(V')$ (the scheme-theoretic image). Note that $V,V'$ and $f(V')$ are integral. As $f$ is proper, $f(V')$ is a closed subscheme of $S$, hence $V\subset f(V')$. As $V_K$ is open in $f(V')$, the equality $f(V')=(f(V')\setminus V_K)\cup V$ implies $f(V')=V$ by the irreducibility of $f(V')$.
\eprf

Now we define the cycle classes on the smaller closed strata. We will see in \S\ref{wssact} that they are essentially equal to the cycle classes $\ol{\Gamma}^{(p)}$ in Proposition 2.20 of \cite{SaT} (for $p=m-1$; we apologize for our numbering of strata which is shifted by 1 from \cite{SaT}). 

\bde \label{cycleclasses}
For two subsets $I,J\subset \Delta$, let $Y_{I,J}:=Y_I\times_kY_J\subset Y\times_k Y$, which is proper smooth of dimension $2n-(|I|+|J|)$ over $k$ if not empty. When $|I|=|J|=m$ and $I=\{i_1,...,i_m\},\ J=\{j_1,...,j_m\}$ both with increasing order, define $D_{I,J}:=D_{i_1,j_1}\cap D_{i_2,j_2}\cap \cdots \cap D_{i_m,j_m}\subset Y'$, which is proper smooth of dimension $2n-1-m$ over $k$ if not empty. As each $D_{i,j}$ is a Cartier divisor on $X'$, we can define
\[\Gamma'_{I,J} := D_{I,J}\cdot \Gamma' := D_{i_1,j_1}\cdot D_{i_2,j_2}\cdots D_{i_m,j_m}\cdot \Gamma'\in A_{n-m}(D_{I,J}\cap |\Gamma'|).\]
Also, as $f(D_{I,J}\cap |\Gamma'|)\subset Y_{I,J}\cap |\Gamma|$ because $f(D_{i,j})=Y_{i,j}$ and $f(|\Gamma'|)=|\Gamma|$, we have:
\[f_*=(f|_{D_{I,J}\cap |\Gamma'|})_*:A_{n-m}(D_{I,J}\cap |\Gamma'|)\lra A_{n-m}(Y_{I,J}\cap |\Gamma|).\]
Finally, let $\Gamma_{I,J}:=f_*(\Gamma'_{I,J})\in A_{n-m}(Y_{I,J}\cap |\Gamma|)$. 

For $1\leq m\leq n$, we set $Y^{(m)}:=\coprod_{|I|=m}Y_I$, and define 
\[ \Gamma^{(m)}:=\sum_{|I|=|J|=m}\Gamma_{I,J}\in A_{n-m}(Y^{(m)}\times_k Y^{(m)}). \] 
\ede

\subsection{Action on the weight spectral sequence} \label{wssact}

Let $\ell$ be a prime $\neq \chara k$, and recall the weight spectral sequence (\cite{RZ},\cite{SaT}) for the $\ell$-adic \'etale cohomology of $X$:
\[E_1^{i,j}:=\bigoplus_{s\geq \max(0,-i)}H^{j-2s}\bigl(Y^{(i+2s+1)}_{\ol{k}}, \Q_\ell(-s)\bigr) \Lra H^{i+j}(X_{\ol{K}}, \Q_\ell).\]
Our interest in the cycle classes defined in Definition \ref{cycleclasses} stems from the following:

\bpr \label{wssaction}
For $1\leq m\leq n$, the action $\oplus [\Gamma^{(i+2s+1)}]^*$ on $E_1^{i,j}$ is compatible with the action $[\Gamma_K]^*$ on $H^{i+j}(X_{\ol{K}}, \Q_\ell)$, i.e.\ we have a map of weight spectral sequences:
\[\xymatrix{
E_1^{i,j}=\bigoplus_{s\geq \max(0,-i)}H^{j-2s}\bigl(Y^{(i+2s+1)}_{\ol{k}},\ \Q_\ell(-s)\bigr) \ar[r]\ar[d]_{\oplus [\Gamma^{(i+2s+1)}]^*} &  H^{i+j}(X_{\ol{K}},\ \Q_\ell) \ar[d]^{[\Gamma_K]^*} \\
E_1^{i,j}:=\bigoplus_{s\geq \max(0,-i)}H^{j-2s}\bigl(Y^{(i+2s+1)}_{\ol{k}},\ \Q_\ell(-s)\bigr) \ar[r] &  H^{i+j}(X_{\ol{K}},\ \Q_\ell)
}\]
\epr

\bprf
This is Proposition 2.20 of \cite{SaT}, except that we need to show that our $[\Gamma^{(i+2s+1)}]^*$ is equal to the $[\ol{\Gamma}^{(p+2s)}]^*/k!$ there. As the cycle class $\ol{\Gamma}^{(p)}$ is defined as the push-forward of $\Gamma^{(p)}$ defined in Lemma 2.17 of \cite{SaT} (restriction $|_{Y_1^{\prime \prime (p)}}$ is simply choosing the relevant connected components), our claim boils down to the following statement: in the Lemma 2.17 of \cite{SaT}, we have a cycle class $c_\ell$ satisfying the desired property without the denominator $k!$, i.e.\ $c_\ell|_{X_K}$ (resp.\ $c_\ell|_{Y^{(p)}}$) is the cycle class of $\Gamma_K$ (resp.\ $\Gamma^{(p)}$), by defining $\Gamma^{(p)}$ as in our Definition \ref{cycleclasses}, i.e.\ $\Gamma^{(p)}:=\sum_{|I|=p+1}(Y_I\cdot \Gamma)$, where $Y_I\cdot \Gamma$ is a cycle class in $Y_I$ defined as $D_{i_1}\cdot D_{i_2}\cdots D_{i_{p+1}}\cdot \Gamma$ if $I=\{i_1,\ldots,i_{p+1}\}$ (this is independent of the ordering by \cite{Fu2}, Corollary 2.4.2). This was proved in \cite{Mi}, \S3.2, using the fact that $\Gamma$ is a flat cycle and then invoking a general theory over a field. 

Here we include an elementary proof using successive intersections of cycles with divisors (\S\ref{intersectdiv}), addressing the Remark 2.18 of \cite{SaT} in this setting. %we prove that the cycle classes defined by Chern characters are compatible with the restrictions $X_K\ra X$ and $Y^{(p)}\ra X$. 
We introduce some notation following \cite{SGA6} XIV, paragraphs 4--6. Let $X$ be a noetherian scheme as in the beginning of \S\ref{review}, and fix $k\in \N$. Let $K^\cdot(X)$ (resp.\ $K_\cdot(X)$) be the Grothendieck group of finite rank locally free (resp.\ coherent) $\mO_X$-modules, and denote its $k$-th graded piece with respect to the $\lambda$-filtration (resp.\ topological filtration, i.e.\ ${\rm Fil}^k_{\rm top}$ consists of classes of sheaves whose support have codimension at least $k$) by ${\rm Gr}^k(X)$ (resp.\ ${\rm Gr}^k_{\rm top}(X)$). Note that the tensor products over $\mO_X$ makes $K^\cdot(X)$ (resp.\ $K_\cdot(X)$) into a ring (resp.\ $K^\cdot(X)$-module). This is compatible with the filtrations (\cite{SGA6} X, 1.3.2); in particular, if $L$ is an invertible sheaf on $X$, then $1-[L]\in {\rm Fil}^1$ and hence tensoring by $L$ acts trivially on ${\rm Gr}^k_{\rm top}(X)$. 

Let $A^k(X)$ be the cycle class group of codimension $k$. By defining $\phi[V]:=[\mO_V]$ for an integral closed subscheme $V$ of $X$ with $[V]\in A^k(X)$, we have a surjective homomorphism $\phi:A^k(X)\ra {\rm Gr}^k_{\rm top}(X)$ (\cite{Fu2}, Example 15.1.5); if $\dim X=n$ and $\mF$ is a coherent $\mO_X$-module whose support has codimension at least $k$, then $\phi$ sends $Z_{n-k}(\mF)$ of \S\ref{sscycles} to $[\mF]$. In particular, (i): if $V$ is a closed subscheme of $X$ with codimension at least $k$, then $\phi(Z_{n-k}(\mO_V))=[\mO_V]$ in ${\rm Gr}^k_{\rm top}(X)$, (ii): if $D$ is a Cartier divisor on $X$, then $\phi[D]=[\mO_X]-[\mO(-D)]$ in ${\rm Gr}^1_{\rm top}(X)$. In fact, if $D$ is effective then $\phi[D]=\phi(Z_{n-1}(\mO_D))=[\mO_D]=[\mO_X]-[\mO(-D)]$. In general, let $D=D_1-D_2$ with $D_1,D_2$ effective. As tensoring by $\mO(-D)$ is trivial on ${\rm Gr}^1_{\rm top}(X)$, we have $[\mO_X]-[\mO(-D_2)]=[\mO(-D)]-[\mO(-D_1)]$. Thus $\phi[D]=\phi[D_1]-\phi[D_2]=([\mO_X]-[\mO(-D_1)])-([\mO(-D)]-[\mO(-D_1)])=[\mO_X]-[\mO(-D)]$.

Now assume $X$ to be regular. The natural homomorphism $K^\cdot(X)\ra K_\cdot(X)$ is an isomorphism with the inverse $[\mF]\mapsto [\mE_\cdot]$, where $\mE_\cdot$ is a locally free resolution of $\mF$ and $[\mE_\cdot]:=\sum_i(-1)^i[\mE_i]$ (see \cite{Fu2}, B.8.3). There is a natural morphism ${\rm Gr}^k(X)\ra {\rm Gr}^k_{\rm top}(X)$ which is an isomorphism modulo torsion (\cite{SGA6} VII, 4.11). From ${\rm Gr}^k(X)$, we have the Chern character ${\rm ch}_k:{\rm Gr}^k(X)\ra H^{2k}(X, \Q_\ell(k))$ (\cite{SGA6} XIV, 5.1), defined via Chern classes of vector bundles. We define the cycle class map as the composite map:
\begin{equation} \label{fact}
\xymatrix{
{\rm cl}: A^k(X) \ar[r]^-\phi & {\rm Gr}^k_{\rm top}(X)/_{\rm (tors)} & {\rm Gr}^k(X)/_{\rm (tors)} \ar[l]_-\cong \ar[r]\ar[r]^-{{\rm ch}_k} & H^{2k}(X,\Q_\ell(k)).
}
\end{equation}
This is equal to the usual cycle class map when $X$ is over a field. It is seen as follows. Both maps factor through $A^k(X)\otimes \Q$. Take any integral closed subscheme $V$ of codimension $k$ in $X$. As we have $\phi[V]=[\mO_V]$ by the above remark (i), we have ${\rm cl}[V]={\rm ch}_k[\mO_V]$. On the other hand, by \cite{Fu2}, Example 15.2.16(a) (essentially the Riemann-Roch theorem for the closed immersions), the Chern character into cycle class groups ${\rm ch}_k: {\rm Gr}^k_{\rm top}(X)\ra A^k(X)\otimes \Q$ satisfies ${\rm ch}_k[\mO_V]=[V]$, and the usual cycle map sends Chern classes to Chern classes (\cite{SGA5} VII, 3.9), hence $[V]$ to ${\rm ch}_k[\mO_V]$. See also \cite{SGA6} XIV, 4.

Back to the proof: we can assume $\Gamma=[V]$ for an integral closed subscheme $V$ of codimension $k$ in $X$, and define $c_\ell:={\rm cl}(\Gamma)$. We need to show $c_\ell|_{X_K}={\rm cl}(\Gamma_K)$ and $c_\ell|_{Y^{(p)}}={\rm cl}(\Gamma^{(p)})$; as $X_K$ and $Y^{(p)}$ are regular and (\ref{fact}) is valid, it suffices to show that (\ref{fact}) is compatible with the restrictions. As $X_K\ra X$ is an open immersion and hence flat, it is clearly compatible with (\ref{fact}). For $Y^{(p)}$, first of all as $Y^{(p)}=\coprod_{|I|=p+1} Y_I$ and $Y_I\cdot \Gamma$ is defined as a successive intersection by the $D_i$, and also at each step we have a cycle on a regular scheme, it suffices to show that (\ref{fact}) is compatible with the restriction to smooth divisors, i.e.\ the following is commutative for a regular $X$ and an effective Cartier divisor $D$ on $X$ which corresponds to a regular closed subscheme:
\[\xymatrix{
A^k(X) \ar[d]_{D\cdot} \ar[r]^-{\rm cl} & H^{2k}(X,\Q_\ell(k)) \ar[d]\\
A^k(D) \ar[r]^-{\rm cl} & H^{2k}(D,\Q_\ell(k))
}\ \ \ \ 
\xymatrix{
[V] \ar@{|->}[d] \ar@{|->}[r] & c_\ell \ar@{|->}[d]\\
D\cdot [V] \ar@{|->}[r] & c_\ell|_D
}
\]
As the Chern character ${\rm ch}_k$ commutes with arbitrary pull-backs (because the Chern classes do; \cite{SGA5} VII, 3.4(ii)), it suffices to show the commutativity of the following:
\[\xymatrix{
A^k(X) \ar[r]^-\phi \ar[d]_{D\cdot} & {\rm Gr}^k_{\rm top}(X) \ar[d] & {\rm Gr}^k(X) \ar[l] \ar[d] \\
A^k(D) \ar[r]^-\phi & {\rm Gr}^k_{\rm top}(D) & {\rm Gr}^k(D) \ar[l]
}\ \ \ \ 
\xymatrix{
[V] \ar@{|->}[d] \ar@{|->}[r] & [\mO_V] \ar@{|->}[d] & [\mE_\cdot] \ar@{|->}[l]  \ar@{|->}[d] \\
D\cdot [V] \ar@{|->}[r]^-{?} & \mF(D,V) & [\mO_D\otimes \mE_\cdot] \ar@{|->}[l]
}
\]
where we set $\mF(D,V):=\bigl[\mO_D\otimes\mO_V\bigr]-\bigl[\ul{\rm Tor}_1(\mO_D,\mO_V)\bigr]$, so that the right square commutes; this is because the cohomology groups of $[\mO_D\otimes \mE_\cdot]$ are $\ul{\rm Tor}_i(\mO_D,\mO_V)$ and hence $[\mO_D\otimes \mE_\cdot]=\sum_i(-1)^i[\ul{\rm Tor}_i(\mO_D,\mO_V)]$ in ${\rm Gr}^k_{\rm top}(D)$, and Lemma \ref{tor}(i). Thus we need to show the commutativity of the left square, i.e.\ the equality $\phi(D\cdot [V])=\mF(D,V)$. As $D\cdot [V]$ and $\mF(D,V)$ both have supports in $D\cap V$ and the maps $\phi$ commute with push-forward by proper morphisms (\cite{Fu2}, Example 15.1.5), in particular the closed immersion $D\cap V\ra D$, it suffices to prove $\phi(D\cdot [V])=\mF(D,V)$ for $\phi:A^j\bigl(D\cap V\bigr)\ra {\rm Gr}^j_{\rm top}\bigl(D\cap V\bigr)$.

We will prove this for any effective Cartier divisor $D$ and integral $V$; note that $D\cdot V=\iota^*D$, where $\iota:V\ra X$ is the closed immersion. We consider the two cases (1) and (2) in \S\ref{intersectdiv} separately. In the case (1) where $V\subset |D|$, we have $D\cap V=V$ and $j=1$, hence we are in ${\rm Gr}^1_{\rm top}(V)$. As $[\iota^*D]$ is a divisor class on $V$ with $\mO(-\iota^*D)=\mO(-D)\otimes \mO_V$, by the remark (ii) before (\ref{fact}) we have $\phi[\iota^*D]=[\mO_V]-[\mO(-\iota^*D)]=[\mO_V]-[\mO(-D)\otimes \mO_V]$. This is equal to $\mF(D,V)$ by the long exact sequence (\ref{torseq}). In the case (2) where $V$ is not contained in $|D|$, by Lemma \ref{tor}(ii) we see $\dim(D\cap V)=i-1$ and $j=0$, and also by Lemma \ref{tor}(ii) we have $\mF(D,V)=[\mO_D\otimes \mO_V]=[\mO_{D\cap V}]$ and $[\iota^*D]=Z_{i-1}(\mO_{D\cap V})$. Thus $\phi[\iota^*D]=\mF(D,V)$ by the remark (i) before (\ref{fact}).
\eprf

\brem
As the proof is essentially a reworking of \cite{SaT}, Lemma 2.17, the above Proposition holds in the setting of \cite{SaT}, Prop.\ 2.20, i.e.\ correspondences of general codimension on $X\times X'$. As the denominator $k!$ is removed from \cite{SaT}, Lemma 2.17/Prop.\ 2.20 (this was also remarked in \cite{Mi}, \S3.2), the proof of $\ell$-independence in \cite{SaT}, \S3 can be somewhat simplified. 
\erem

\subsection{A formula relating the cycle classes on the strata}

We return to the notation of Definition \ref{cycleclasses}.

\bpr \label{iteratedproj}
When $|I|=|J|+1=m$ and $I=\{i_1,...,i_m\},\ J=\{j_1,\ldots,j_{m-1}\}$ are in increasing order, we have the equality
\[ \sum_{h=1}^m(-1)^h\,Y_{I,J}\cdot \Gamma_{I\setminus \{i_h\},J} = \sum_{j\in \Delta\setminus J} (-1)^{h(j)}\,\Gamma_{I,J\cup\{j\}} \]
in $A_{n-m}(Y_{I,J}\cap |\Gamma|)$, where $1\leq h(j)\leq m$ is determined by $j_{h(j)-1}<j<j_{h(j)}$ (set $j_m:=\infty$), i.e.\ the position of $j$ in $J\cup \{j\}$ when the elements are ordered increasingly.
\epr

\bprf
First we extend the definition of $D_{I,J\amalg \{j\}}, Y_{I,J\amalg\{j\}}$ and $\Gamma_{I,J\amalg\{j\}}$ to the case $j\in J$: if $J\amalg\{j\}=\{j_1,...,j_m\}$ with increasing order (with double $j$), define $D_{I,J\amalg \{j\}}:=D_{i_1,j_1}\cap D_{i_2,j_2}\cap \cdots \cap D_{i_m,j_m}$. Then the proposition follows from the following lemma by taking the alternating sum with respect to $h$, as each term with $j\in J$ cancel out by appearing in two successive values of $h$.
\eprf

\blem
For $1\leq h\leq m$, we have the equality in $A_{n-m}(Y_{I,J}\cap |\Gamma|)$:
\[ Y_{I,J}\cdot \Gamma_{I\setminus\{i_h\},J}=\sum_{j_{h-1}\leq j\leq j_h}\Gamma_{I,J\amalg\{j\}}. \]
\elem

\bprf
Write $f=f|_{D_{I\setminus\{i_h\},J}}: D_{I\setminus\{i_h\},J}\lra Y_{I\setminus\{i_h\},J}$. As $\Gamma_{I\setminus\{i_h\},J}=f_*(D_{I\setminus\{i_h\},J}\cdot \Gamma')$ by definition, considering $Y_{I,J}$ as a Cartier divisor of $Y_{I\setminus\{i_h\},J}$, the projection formula (Proposition \ref{projformula}(vi)) reads:
\[Y_{I,J}\cdot \Gamma_{I\setminus\{i_h\},J} = g_*(f^*Y_{I,J}\cdot \Gamma')\ \ \text{ in } A_{n-m}(Y_{I,J}\cap |\Gamma|),\]
where $g:=f|_{f^{-1}(Y_{I,J})\cap |\Gamma|}$. Consider the commutative diagram:
\[\xymatrix{
D_{I\setminus\{i_h\},J} \ar[d]^-f \ar[r]^-{i'} & X' \ar[d]^-f\\
Y_{I\setminus\{i_h\},J} \ar[r]^-i & X\times_\mO X 
}\]
Using Lemma \ref{cartdivsum}, we have the equality of Cartier divisors:
\[ f^*Y_{I,J}=f^*i^*\pr_1^*Y_{i_h}=i^{\prime *}f^*\pr_1^*Y_{i_h}=i^{\prime *}\sum_{j\in \Delta}D_{i_h,j}, \]
but as $D_{I\setminus\{i_h\},J}\cap D_{i_h,j}=\emptyset$ unless $j_{h(j)-1}\leq j\leq j_{h(j)}$ by Proposition \ref{basicblowup}(iii), we have 
\[ i^{\prime *}\sum_{j\in \Delta}D_{i_h,j} \ \ = \sum_{j_{h(j)-1}\leq j\leq j_{h(j)}}\!\!\!D_{i_h,j}\cdot D_{I\setminus\{i_h\},J} \ \ = \sum_{j_{h(j)-1}\leq j\leq j_{h(j)}}\!\!\!D_{I,J\amalg\{j\}}. \]
Therefore combining with the projection formula we obtain the lemma.
\eprf

\brem
Clearly the formula depends on the ordering of $\Delta$, i.e.\ the components $Y_i$, but so do the blow up $X'$ and the Definition \ref{cycleclasses} (and perhaps the weight spectral sequence itself). The author does not know how canonical these constructions can be made.
\erem

\section{Proof of the main results}

\subsection{Characterization of correspondences} \label{charcorres}

Now we would like to use Proposition \ref{iteratedproj} to compute $\Gamma^{(m)}=\sum_{|I|=|J|=m} \Gamma_{I,J}$ inductively on $m$. In order to do this, we need to identify $\Gamma_{I,J}$ not only as cycle classes, but as cycles.

We begin with the case $|I|=|J|=1$. For $i,j\in \Delta$, let $\Gamma'_{i,j}:=\Gamma'_{\{i\},\{j\}}=D_{i,j}\cdot \Gamma'$. As $\Gamma'$ is the closure of $\Gamma_K$, it is not contained in $D_{i,j}$, therefore $\Gamma'_{i,j}$ is well-defined as an element of $Z_{n-1}(D_{i,j}\cap |\Gamma'|)$. Correspondingly, $\Gamma_{i,j}=f_*(\Gamma'_{i,j})$ is a well-defined element of $Z_{n-1}(Y_{i,j}\cap |\Gamma|)$. 

To proceed further, we will make the following assumption on $\Gamma_K$. Recall that $|\Gamma|$ is a reduced closed subscheme of $X\times_\mO X$. We assume:
\[(*)\ \ \ \ \ 
\text{ Two projections } \pr_1,\pr_2:|\Gamma|\lra X \text{ are finite morphisms.}\]
The immediate consequence of this assumption is that
\[\dim(Y_{I,J}\cap |\Gamma|)\leq \min\{n-|I|,\ n-|J|\},\]
because the projection morphisms $Y_{I,J}\cap |\Gamma|\lra Y_I$ and $Y_{I,J}\cap |\Gamma|\lra Y_J$ are both finite and $\dim Y_I=n-|I|$. In particular, when $|I|=|J|=m$, we have $\Gamma_{I,J}\in A_{n-m}(Y_{I,J}\cap |\Gamma|)=Z_{n-m}(Y_{I,J}\cap |\Gamma|)$. Moreover, the equality in Proposition \ref{iteratedproj} is an equality in $Z_{n-m}(Y_{I,J}\cap |\Gamma|)$, and as the cycle on the RHS belongs to its subgroup
\[Z_{n-m}\Bigl(\bigcup_{j\in \Delta\setminus J}(Y_{I,J\cup\{j\}}\cap |\Gamma|)\Bigr)=\bigoplus_{j\in \Delta\setminus J}Z_{n-m}(Y_{I,J\cup\{j\}}\cap |\Gamma|).\]
(This equality follows from $\dim (Y_{I,J\cup\{j\}}\cap Y_{I,J\cup\{j'\}}) = \dim (Y_{I,J\cup\{j,j'\}}) < n-m$.) Therefore so does the cycle on the LHS. This means that if we know the set of cycles $\Gamma_{I,J}$ when $|I|=|J|=m-1$, then computing the LHS yields the set of cycles $\Gamma_{I,J}$ when $|I|=|J|=m$. This reduces the determination of $\Gamma_{I,J}$ to the case $|I|=|J|=1$. 

Now recall that $f:X'\ra X\times_\mO X$ is an isomorphism outside $Z=X_{\rm sing} \times_k X_{\rm sing}$, and $D_{i,j}$ is the closure of $Y_{i,j}\setminus Z$ in $X'$. Therefore $f$ is an isomorphism on
\[ (X\times_\mO  X)_{\rm sm}:=X_{\rm sm}\times_{\mO}X_{\rm sm}=(X\times_{\mO}X)\setminus Z',\ \ \text{ where } Z':=\!\!\bigcup_{\max(|I|,|J|)\geq 2}\!\!\!Y_{I,J}.\]
Therefore if we write $Y_i^0:=Y_i\cap X_{\rm sm}$, then $Y_{i,j}^0:=Y_i^0\times_kY_j^0=Y_{i,j}\setminus Z'$ is a Cartier divisor of $(X\times_{\mO} X)_{\rm sm}$ (which is also clear from $(X\times_{\mO} X)_{\rm sm}\times_\mO k=\coprod_{(i,j)\in \Delta\times \Delta}Y_{i,j}^0$). As $Y_{i,j}^0$ is open dense in $Y_{i,j}$, the $D_{i,j}$ is also the closure of $Y_{i,j}^0$ in $X'$.

As $f$ is an isomorphism on $(X\times_{\mO} X)_{\rm sm}$, we can identify $f^{-1}(Y_{i,j}^0)\cap |\Gamma'|$ with $Y_{i,j}^0\cap |\Gamma|$ by $f$, and by Proposition \ref{projformula}(vii):
\[ \Gamma_{i,j}^0:=\Gamma_{i,j}|_{Y_{i,j}^0\cap |\Gamma|}=\Gamma'_{i,j}|_{f^{-1}(Y_{i,j}^0)\cap |\Gamma'|}=Y_{i,j}^0\cdot \bigl(\Gamma|_{(X\times_{\mO}X)_{\rm sm}}\bigr)\in Z_{n-1}(Y_{i,j}^0\cap |\Gamma|).\]
Now by our assumption $(*)$, we have $\dim (Z'\cap |\Gamma|)\leq n-2$, therefore $Z_{n-1}(Y_{i,j}\cap |\Gamma|)\liso Z_{n-1}(Y_{i,j}^0\cap |\Gamma|)$ by Lemma \ref{opencycles}, and $\Gamma_{i,j}$ must equal the closure of $\Gamma^0_{i,j}$ in $Y_{i,j}\cap |\Gamma|$. Summarizing, we obtain a characterization of $\Gamma_{I,J}$ which does not involve the blow-up $f:X'\lra X\times_{\mO} X$:

\bpr \label{characterization}
Under $(*)$, the collection $\{\Gamma_{I,J}\}$ of cycles $\Gamma_{I,J}\in Z_{n-|I|}(Y_{I,J}\cap |\Gamma|)$ for all pairs $(I,J)$ with $|I|=|J|$ is characterized by the following two properties:
\benu
\item $\Gamma_{i,j}$ is the closure of $\Gamma_{i,j}^0:=Y_{i,j}^0\cdot \Gamma|_{(X\times_{\mO}X)_{\rm sm}}$ in $Y_{i,j}$. 
\item The equality in Proposition \ref{iteratedproj} for all pairs $(I,J)$ with $|I|=|J|+1$.
\eenu
Moreover, the equality in Proposition \ref{iteratedproj} is the equality of cycles.
\epr

This together with Proposition \ref{wssaction} proves the first part of the Theorem in the introduction.

\subsection{Composing the correspondences} \label{composesec}

For two correspondences $\Gamma_1,\Gamma_2$ on a noetherian scheme $X$ which is proper smooth over $S=\Spec K,\ \Spec k\text{ or }\Spec \mO$, their {\em composition} is defined as (see \cite{Fu2}, \S16.1):
\[ \Gamma_2\circ \Gamma_1:=p_{13*}(p_{23}^*\Gamma_2\cdot p_{12}^*\Gamma_1)\in A\bigl(p_{13}((X\times |\Gamma_2|)\cap (|\Gamma_1|\times X))\bigr), \]
where $p_{ij}:X\times_S X \times_S X \lra X\times_S X$ denotes the projection onto $i$-th and $j$-th components and $\cdot$ is the intersection product of cycle classes (\cite{Fu2}, Chapters 6 and 8 over a field, and \S20.2 for smooth schemes over a regular 1-dimensional base). This composition induces the composition of their actions on the cohomology groups (\S\ref{algcor}).
%\[ [\Gamma_2\circ \Gamma_1]^* = [\Gamma_1]^* \circ [\Gamma_2]^*. \]

Now we return to our situation of strictly semistable scheme $X$ over $\mO$ and let $\Gamma_{1,K},\Gamma_{2,K}$ be correspondences on $X_K$. If $\Gamma_K:=\Gamma_{2,K}\circ \Gamma_{1,K}$, then:
\[ [\Gamma_K]^* = [\Gamma_{1,K}]^* \circ [\Gamma_{2,K}]^* \ \text{ on }\ H^i(X_{\ol{K}}, \Q_\ell).\]
We denote the closures of $\Gamma_{1,K}, \Gamma_{2,K}$ and $\Gamma_K$ in $X\times_\mO X$ respectively by $\Gamma_1,\Gamma_2$ and $\Gamma$. Then one might conjecture that the family of cycle classes $\Gamma^{(m)}$ for $1\leq m\leq n$ in Definition \ref{cycleclasses} should be compatible with compositions of correspondences, i.e.\ $\Gamma^{(m)}=\Gamma_2^{(m)}\circ \Gamma_1^{(m)}$. The author does not have a proof in general, but under the assumption $(*)$, we can prove its validity using Proposition \ref{characterization}. First note the following:

\blem \label{composefin}
Assume that $\Gamma_{1,K},\Gamma_{2,K}$ both satisfy $(*)$. Then:
\benu
\item $\Gamma_K:=\Gamma_{2,K}\circ \Gamma_{1,K}$ satisfies $(*)$. Also, it is well-defined as an $(n-1)$-dimensional cycle on $X_K\times X_K$.
\item the definition of $\Gamma_{2,K}\circ \Gamma_{1,K}$ extends to $X_{\rm sm}$ (smooth but {\em not} proper over $\mO$). In particular on $Y^0:=\coprod_{i,j}Y_{i,j}^0$ we have $\Gamma|_{Y^0\cap |\Gamma|}=\Gamma_2|_{Y^0}\circ \Gamma_1|_{Y^0}$.
%$\Gamma|_{(X\times_{\mO}X)_{\rm sm}\cap |\Gamma|}=\Gamma_2|_{(X\times_{\mO}X)_{\rm sm}} \circ \Gamma_1|_{(X\times_{\mO}X)_{\rm sm}}$.
\eenu
\elem

\bprf
(i): We put the indices $|\Gamma_1|\subset X_1\times X_2$ and $|\Gamma_2|\subset X_2\times X_3$ to distinguish the different projections (products are over $\mO$), and let $Y:=(X_1\times |\Gamma_2|)\cap (|\Gamma_1|\times X_3)$, a closed subscheme of $X_1\times X_2\times X_3$. We show that the projections $Y\ra X_1,\ Y\ra X_3$ are finite.
\[\xymatrix{
Y \ar[r]\ar[d]_-{\text{cl.imm.}} & |\Gamma_1| \ar[dr]^-{\text{fin.}}\ar[d]_-{\text{cl.imm.}} & \\
X_1\times |\Gamma_2| \ar[r]^-{\text{fin.}} & X_1\times X_2 \ar[r] & X_1
}\]
As $\Gamma_2$ satisfies $(*)$, we know that $X_1\times |\Gamma_2|\ra X_1\times X_2$ is finite. As $Y\subset X_1\times |\Gamma_2|$, the morphism $Y\ra X_1\times X_2$ is finite. It factors through $Y\ra |\Gamma_1|$, hence this morphism is also finite. Hence by $(*)$ for $\Gamma_1$ we conclude that $Y\ra X_1$ is finite. Similarly $Y\ra X_3$ is finite. Thus two projections $p_{13}(Y)\ra X_1,\ p_{13}(Y)\ra X_3$ are both finite. As $|\Gamma_K|$ is a closed subscheme of $p_{13}(Y\times_\mO K)$, its closure $|\Gamma|$ is a closed subscheme of $p_{13}(Y)$, hence the projections $|\Gamma|\ra X_1,\ |\Gamma|\ra X_3$ are both finite. (ii): We proved that $p_{13}:Y\lra p_{13}(Y)$ is finite, thus $p_{13*}$ makes sense. As $Y^0=(X\times_\mO X)_{\rm sm}\times_\mO k$, the equality follows from the compatibility of the intersection product with the specialization (\cite{Fu2}, \S20.2 and \S20.3).
\eprf

\bpr \label{compose}
Assume that $\Gamma_{1,K},\Gamma_{2,K}$ both satisfy $(*)$, and let $\Gamma$ be the closure of $\Gamma_{2,K} \circ \Gamma_{1,K}$ in $X\times_\mO X$. Then we have $\Gamma^{(m)}=\Gamma_2^{(m)}\circ \Gamma_1^{(m)}$ for $1\leq m\leq n$. In particular, we have $[\Gamma^{(m)}]^* = [\Gamma_1^{(m)}]^* \circ [\Gamma_2^{(m)}]^*$ on the cohomology groups of $Y^{(m)}$.
\epr

\bprf
In this proof only we use $K$ and $k$ for a subset and an element of $\Delta$, and write:
\[ \Gamma_1=\sum_{|I|=|J|=m}\Gamma^1_{I,J},\ \ \ \ \Gamma_2=\sum_{|J|=|K|=m}\Gamma^2_{J,K},\ \ \ \ \Gamma=\sum_{|I|=|K|=m}\Gamma_{I,K}.\] 
By Lemma \ref{composefin}(i) and Proposition \ref{characterization}, it suffices to check that $\Gamma_2^{(m)}\circ \Gamma_1^{(m)}$ satisfies the two conditions for $(\Gamma_2\circ \Gamma_1)^{(m)}$ in Proposition \ref{characterization}. We know the condition for $m=1$ by Lemma \ref{composefin}(ii). 
%, setting $Y^0:=\coprod_{i,j}Y_{i,j}^0$, we need to check that the cycles $\Gamma|_{Y^0\cap |\Gamma|}$ are the same as $\Gamma_2|_{Y^0}\circ \Gamma_1|_{Y^0}$, which follows from Lemma \ref{composefin}(ii) and
By induction on $m$, assume that we know $\Gamma^{(m-1)}=(\Gamma_2\circ \Gamma_1)^{(m-1)}$, i.e.:\ 
\[ \Gamma_{I,K}=\sum_{|J|=m-1}\Gamma^2_{J,K}\circ \Gamma^1_{I,J}\ \text{ for }\ |I|=|K|=m-1.\]
We use the notation $Y_{I,J,K}:=Y_I\times Y_J\times Y_K$. For $|I|=m$ and $|K|=m-1$, we compute as follows (we are on proper smooth schemes over $k$ now):
\begin{align*}
\sum_{h=1}^m(-1)^hY_{I,K}\cdot \Gamma_{I\setminus\{i_h\},K} 
&= \sum_{h=1}^m(-1)^h\sum_{|J|=m-1}Y_{I,K}\cdot (\Gamma^2_{J,K}\circ \Gamma^1_{I\setminus\{i_h\},J})\\
&= \sum_{h=1}^m(-1)^h\sum_{|J|=m-1}Y_{I,K}\cdot p_{13*} \bigl((Y_{I\setminus\{i_h\}}\times \Gamma^2_{J,K})\cdot (\Gamma^1_{I\setminus\{i_h\},J} \times Y_K) \bigr)\\
&= \sum_{h=1}^m(-1)^h p_{13*} \sum_{|J|=m-1} Y_{I,J,K}\cdot (Y_{I\setminus\{i_h\}}\times \Gamma^2_{J,K})\cdot (\Gamma^1_{I\setminus\{i_h\},J} \times Y_K),
%&= p_{13*} \sum_{h=1}^m(-1)^h \sum_{|J|=m-1} (Y_I\times \Gamma^2_{J,K})\cdot (Y_{I,J}\cdot \Gamma^1_{I\setminus\{i_h\},J} \times Y_K)
\end{align*}
where we used the projection formula. When $W$ is a smooth hypersurface not containing $V_1,V_2$ we have $[W]\cdot[V_1]\cdot[V_2]=[W\cdot V_1]\cdot [W\cdot V_2]$ (see \cite{Fu2} Example 8.1.10), therefore
\[ \text{ RHS } = p_{13*}\sum_{|J|=m-1} (Y_I\times \Gamma^2_{J,K})\cdot \bigl( \sum_{h=1}^m(-1)^h Y_{I,J}\cdot \Gamma^1_{I\setminus\{i_h\},J} \times Y_K \bigr). \]
%&= p_{13*}\sum_{|J|=m-1} (Y_I\times \Gamma^2_{J,K})\cdot \bigl( \sum_{j\in \Delta\setminus J}(-1)^{h(j)}\Gamma^1_{I,J\cup\{j\}} \times Y_K \bigr)\\
Now using Proposition \ref{characterization} for $\Gamma_1$ and reindexing: 
\begin{align*}
\text{ RHS } 
&= p_{13*} \sum_{|J|=m-1} \sum_{j\in \Delta\setminus J}(-1)^{h(j)} (Y_I\times \Gamma^2_{J,K})\cdot (\Gamma^1_{I,J\cup\{j\}} \times Y_K)\\
&= p_{13*} \sum_{|J'|=m} \sum_{h=1}^{m+1} (-1)^h (Y_I\times \Gamma^2_{J'\setminus \{j_h\},K})\cdot (\Gamma^1_{I,J'} \times Y_K),
\end{align*}
then working backwards as above and using Proposition \ref{characterization} for $\Gamma_2$ gives:
\begin{align*}
\text{ RHS }
&= p_{13*} \sum_{|J'|=m} \sum_{h=1}^{m+1} (-1)^h Y_{I,J',K}\cdot (Y_I\times \Gamma^2_{J'\setminus \{j_h\},K})\cdot (\Gamma^1_{I,J'} \times Y_K)\\
&= p_{13*} \sum_{|J'|=m} (Y_I\times \sum_{h=1}^{m+1} (-1)^h Y_{J',K}\cdot \Gamma^2_{J'\setminus \{j_h\},K}) \cdot (\Gamma^1_{I,J'} \times Y_K)\\
%&= p_{13*} \sum_{|J'|=m} (Y_I\times \sum_{k\in \Delta\setminus K} (-1)^{h(k)} \Gamma^2_{J',K\cup \{k\}}) \cdot (\Gamma^1_{I,J'} \times Y_K)\\
&= p_{13*} \sum_{|J'|=m} \sum_{k\in \Delta\setminus K} (-1)^{h(k)} (Y_I\times \Gamma^2_{J',K\cup \{k\}}) \cdot (\Gamma^1_{I,J'} \times Y_K) \\
&= \sum_{k\in \Delta\setminus K} (-1)^{h(k)} \sum_{|J'|=m} \Gamma^2_{J',K\cup \{k\}}\circ \Gamma^1_{I,J'},
\end{align*}
which is equal to $\sum_{k\in \Delta\setminus K} (-1)^{h(k)} \Gamma_{I,K\cup \{k\}}$ as desired.
\eprf

This proves the second part of the Theorem in the introduction.

\end{document}